\documentclass[12pt,a4paper]{article}

\usepackage[dutch,english]{babel}
\usepackage[UKenglish]{isodate}
\usepackage{amssymb, amsmath, calligra, mathrsfs, color, fancyhdr, enumitem, amsthm, hyperref, graphicx, picins, array, tikz-cd, subfiles, footnote}
\usepackage[toc,page]{appendix}
\usepackage[a4paper, top=3.5cm, bottom=3cm, left=3cm, right=3cm]{geometry}
\usepackage[all]{xy}
\newtheorem{theorem}{Theorem}
\newtheorem{lemma}[theorem]{Lemma}

\newtheorem{conjecture}[theorem]{Conjecture}
\theoremstyle{definition}
\newtheorem{definition}[theorem]{Definition}
\newtheorem{remark}[theorem]{Remark}

\newtheorem{algorithm}[theorem]{Algorithm}
\newtheorem{notation}[theorem]{Notation}
\renewcommand{\O}{\mathcal{O}}

\newcommand{\commentaar}[1]{}
\newcommand{\C}{\ensuremath{\mathbb{C}}}
\newcommand{\R}{\ensuremath{\mathbb{R}}}
\newcommand{\Q}{\ensuremath{\mathbb{Q}}}
\newcommand{\Z}{\ensuremath{\mathbb{Z}}}

\newcommand{\F}{\ensuremath{\mathbb{F}}}
\newcommand{\A}{\ensuremath{\mathbb{A}}}

\newcommand{\Hom}{\textrm{Hom}}
\DeclareMathOperator{\ShHom}{\mathscr{H}\text{\kern -3pt {\calligra\large om}}\,}
\renewcommand{\O}{\mathcal{O}}

\newcommand{\Pic}{\mathrm{Pic}}

\newcommand{\Spec}{\mathop\mathrm{Spec}}

\newcommand{\tors}{\mathrm{tors}}
\makesavenoteenv{tabular}
\newcounter{nootje}
\renewcommand\check[1]{[*\thenootje]\marginpar{\tiny\begin{minipage}{25mm}\begin{flushleft}\thenootje : #1\end{flushleft}\end{minipage}}\addtocounter{nootje}{1}}
\renewcommand\check[1]{}

\DeclareFontFamily{U}{wncy}{}
\DeclareFontShape{U}{wncy}{m}{n}{<->wncyr10}{}
\DeclareSymbolFont{mcy}{U}{wncy}{m}{n}
\DeclareMathSymbol{\Sh}{\mathord}{mcy}{"58} 

\newcommand{\onlyinsubfile}[1]{#1}

\begin{document}

\title{Numerical verification of the Birch and Swinnerton-Dyer conjecture for  hyperelliptic curves of higher genus over $\Q$ up to squares}
\author{Raymond van Bommel\footnote{Mathematisch Instituut, Universiteit Leiden, P.O.\ box 9512, 2300RA Leiden, The Netherlands}\\ {\footnotesize r.van.bommel@math.leidenuniv.nl}}
\cleanlookdateon
\maketitle

%{\color{red} Idea for extension: find an elliptic curve over a quadratic number field (say $\Q(i)$, $\Q(\zeta_3)$ or $\Q(\sqrt{2})$) that is isogenous to the Jacobian of a genus 2 curve for which we can verify BSD. See for example \href{https://www.impan.pl/pl/wydawnictwa/czasopisma-i-serie-wydawnicze/acta-arithmetica/all/116/3/81940/abelian-surfaces-of-rm-gl-2-type-as-jacobians-of-curves}{this} and \href{https://pdfs.semanticscholar.org/89d2/bd161d5096c6c78942eb5e46e564caf3f9ae.pdf}{this} text.}

\renewcommand{\onlyinsubfile}[1]{}

\textbf{Abstract.} The Birch and Swinnerton-Dyer conjecture has been numerically verified for the Jacobians of 32 modular hyperelliptic curves of genus 2 by Flynn, Lepr\'evost, Schaefer, Stein, Stoll and Wetherell, using modular methods. In the calculation of the real period, there is a slight inaccuracy, which might give problems for curves with non-reduced components in the special fibre of their N\'eron model. In this present paper we explain how the real period can be computed, and how the verification has been extended to many more hyperelliptic curves, some of genus 3, 4 and 5, without using modular methods.


\section{Introduction}			% STATEMENT AND GOAL %

In \cite{BirchSwinnertonDyer}, Birch and Swinnerton-Dyer first stated their famous conjecture, based on computations with elliptic curves. Later, in \cite{TateBourbaki}, Tate generalised the conjecture to abelian varieties of higher dimension.

\begin{conjecture}[\textrm{BSD, \cite[Conj.\ F.4.1.6, p.\ 462]{HindrySilverman}}]
Let $A/\Q$ be an abelian variety of dimension $d$ and algebraic rank $r$. Let $L(A,s)$ be its $L$-function, $A^{\vee}$ its dual, $R_A$ its regulator, $\Sh(A)$ its Tate-Shafarevich group and $P_A$ its period. For each prime $p$, let $c_p$ be the Tamagawa number of $A$ at $p$.
Then $L(A,s)$ has a zero of order $r$ at $s = 1$ and $$\lim_{s \rightarrow 1} {(s-1)^{-r} L(A,s)} = \frac{P_AR_A \cdot |\Sh(A)| \cdot \prod_{p} c_p}{|A(\Q)_{\tors}| \cdot |A^\vee(\Q)_{\tors}|}.$$
\end{conjecture}

\begin{remark}
In Tate's original version, \cite{TateBourbaki}, the period, Tamagawa numbers and discriminant are put in the normalisation of the $L$-function. 
\end{remark}

Tate stated the conjecture for abelian varieties over number fields. However, in \cite{MilneInventiones}, Milne proved that the conjecture is compatible with Weil restriction, so BSD holds for all abelian varieties over all number fields if and only if it holds for all abelian varieties over $\Q$.

Due to work of Kolyvagin (\cite{Kolyvagin}, \cite{KolyvaginSurvey}) and others, a weak version of BSD has been proven for elliptic curves over $\Q$ with analytic rank at most 1. More precisely, we know that in these cases the algebraic rank equals the analytic rank. On the other hand, on the numerical side, in \cite{EmpiricalEvidence} Flynn et al.\ numerically verified BSD for the Jacobians of 32 hyperelliptic curves of genus 2 with small conductor, using modular methods for their calculations.

There is, however, a slight inaccuracy in \cite{EmpiricalEvidence}. In the calculation of the real period, calculations seem to be done inside the sheaf of relative differentials, while they should be done inside the canonical sheaf. For curves whose N\'eron model has non-reduced fibres, this could cause a problem. For the curves considered, it did not seem to invalidate the final results.

The goal of this paper is twofold. On the one hand, we will give a more explicit algorithm to compute the real period, or more specifically, a N\'eron differential, along with the theoretical foundations that are needed for this. On the other hand, we will present how we extended the numerical verification of BSD to many more hyperelliptic curves of genus 2, 3, 4 and 5 without using modular methods. As far as the author is aware, this is the first time BSD has been numerically verified for curves of genus 3, 4 and 5.

We did not compute, however, the order of $\Sh(A)$. Moreover, the verification is only provable up to squares. That is, all terms but $|\Sh(A)|$ are computed, of which some are only provably correct up to squares. Then it is verified that the conjectural order of $\Sh(A)$, as predicted by the conjecture, up to a certain high precision, is a rational square or two times a rational square, in accordance with the criteria described in \cite{PoonenStoll}.

The structure of this article is as follows. First we present our verification results. Then we discuss the computation of the real period and the theoretical background needed. Then in the last part we briefly discuss the computation of the other terms in the BSD formula.

The author wishes to thank his supervisors David Holmes and Fabien Pazuki, Tim Dokchitser, Steffen M\"uller, Carlo Pagano and an anonymous referee for the comments they provided to improve this paper.

%\begin{example}
%Let $H_1$ be the hyperelliptic curve over $\Q$ defined by $$y^2 + (x^3 + x + 1)y = x^6 + 5x^5 + 12x^4 + 12x^3 + 6x^2 - 3x - 4.$$ It is curve 125,B from \cite{EmpiricalEvidence}. Its discriminant is $5^{16}$ and its conductor is 125. We will numerically verify BSD for $J_1$, the Jacobian of $H_1$, and check that the results agree with those from \cite{EmpiricalEvidence}.
%\end{example}

%\begin{example}
%Let $H_2$ be the hyperelliptic curve over $\Q$ defined by $$y^2 = x^5 - 2x^4 - 2x^3 + 4x^2 + x - 1,$$ having discriminant $62720 = 2^8 \cdot 5 \cdot 7^2$ and conductor $7840$ (label $7840.a.62720.1$ from \cite{lmfdb}). We will numerically verify BSD for $J_2$, the Jacobian of $H_2$, which has algebraic rank equal to 1.
%\end{example}

%\begin{example}
%Let $H_3$ be the hyperelliptic curve over $\Q$ defined by $$y^2 = x^7 - x^6 + 3x^5 - x^4 + 2x^3 + x^2 + 1,$$ having discriminant $-1523712 = -2^{14} \cdot 3 \cdot 31$. We will numerically verify BSD for $J_3$, the Jacobian of $H_3$, which has algebraic rank 1.
%\end{example}

\onlyinsubfile{

}

\section{Results}		% NUMERICAL VERIFICATION %

For the Jacobians of the curves listed below, we numerically verified BSD in the following sense. We numerically determined the algebraic and analytic rank, the special value of the $L$-function, the regulator (provably only up to squares), the real period, the Tamagawa numbers, and the size of the torsion subgroup of the Jacobian, assuming some conjectures mentioned below. Then the BSD formula was used to calculate a conjectural order for $\Sh$, and it was verified that it is a rational square (which it should be according to the criteria in \cite{PoonenStoll}).

In practice this meant that the conjectural order for $\Sh$ was less than $10^{-9}$ away from an integer. Moreover, for all but one of the curves of genus 2, this conjectural order was actually equal to 1.000000000.

The conjectural results that we assume to hold for the verification include the analytic continuation, and the correctness of the functional equation of the $L$-function (see \cite[Conj.\ F.4.1.5, p.\ 461]{HindrySilverman}). When we computed the analytic rank, we did this by numerically checking whether the $L$-function and its derivatives up to certain order, vanish at 1. Even though this does not prove that these functions vanish, we do assume this to be true. Moreover, we assume the correctness of Ogg's formula for the computation of the 2-part of the conductor (for more details, see Remark \ref{rmk:Ogg}). In a certain sense, one could say that our verification also provides evidence for these conjectures.

{\bf List of curves:}\\[-0.8cm]
\begin{itemize}
\item All elliptic curves of the form $y^2 = x^3 + ax + b$ with $a,b \in \{-15, \ldots, 15\}$, and compared it with the outcomes of already existing algorithms in \texttt{Magma}.
\item All hyperelliptic curves from \cite{EmpiricalEvidence}, comparing it with the outcomes given in that article.
\item All 300 hyperelliptic curves $C$ of genus 2, of the form $$y^2 = x^5 + ax^4 + bx^3 + cx^2 + dx + e,$$ up to isomorphism, with $a,b,c,d,e \in \{-10, \ldots, 10\}$ and $\Delta(C) \leq 10^5$. About one third of them have rank 1, the rest are of rank 0. They are all contained in the LMFDB, cf.\ \cite{Database2}.
\item All 6 hyperelliptic curves of genus 3, of the form $$y^2 = x^7 + ax^6 + bx^5 + cx^4 + dx^3 + ex^2 + fx + g,$$ with $a,b,c,d,e,f,g \in \{-3, \ldots, 3\}$ and $\Delta(C) \leq 10^6$,
i.e., we checked BSD, up to squares, for
\begin{itemize}
\item $H_1\colon (a,b,c,d,e,f,g) = (1, -3, 2, 2, -3, 0, 1, 0),$
\item $H_2\colon (a,b,c,d,e,f,g) = (1, -2, -1, 2, 2, -1, -1, 0),$
\item $H_3\colon (a,b,c,d,e,f,g) = (1, 0, -3, -2, 2, 3, 1, 0),$
\item $H_4\colon (a,b,c,d,e,f,g) = (1, 0, -1, 0, -2, 3, -1, 0),$
\item $H_5\colon (a,b,c,d,e,f,g) = (1, 1, -2, -2, 1, 2, -1, 0),$
\item $H_6\colon (a,b,c,d,e,f,g) = (1, -3, 2, 0, 1, 0, -1, 0),$
\end{itemize}
and, in order to have an example of rank 1, the curve
\begin{itemize}
\item $H_7\colon (a,b,c,d,e,f,g) = (1,-3,1,3,-2,0,1,0).$
\end{itemize}
As far as we are aware these are the first examples of curves of genus 3 for which BSD has been numerically verified. These were the invariants we found:

\scalebox{0.86}{
\begin{tabular}{|c|c|c|c|c|c|c|c|}\hline
		&$r$	&$\lim_{s \rightarrow 1} \ldots$	&$P_A$	&$R_A$	&$c_p$	&$|A(\Q)_{\mathrm{tors}}|$	&$|\Sh|_{\mathrm{an}}$	\\ \hline
$H_1$	&0		&0.8006061	&51.23879	&1			&$c_2 = c_5 = c_{23} = 1$		&8	&1.000000	\\ \hline
$H_2$	&0		&0.7636550	&48.87392	&1			&$c_2 = c_5 = c_{23} = 1$		&8	&1.000000	\\ \hline
$H_3$	&0		&0.9275079	&59.36050	&1			&$c_2 = c_5 = c_{23} = 1$		&8	&1.000000	\\ \hline
$H_4$	&0		&0.8087909	&51.76262	&1			&$c_2 = c_5 = c_{31} = 1$		&8	&1.000000	\\ \hline
$H_5$	&0		&0.9784790	&62.62265	&1			&$c_2 = c_5 = c_{23} = 1$		&8	&1.000000	\\ \hline
$H_6$	&0		&0.4310775	&55.17793	&1			&$c_2 = 2, c_5 = c_{23} = 1$	&16	&1.000000	\\ \hline
$H_7$	&1		&1.953631	&50.85263	&0.6146799	&$c_2 = c_5 = c_{11} = 1$		&4	&1.000000	\\ \hline
\end{tabular}}

For the torsion and regulator, points were searched up to a certain height on the Jacobian. This maximum search height is considerably smaller than the height given by the various height bounds in the literature. It is possible that the size of the torsion subgroups and the regulator is incorrect, but this would only cause a rational square error factor for the value of $|\Sh|_{\mathrm{an}}$.
\item The curve $$y^2 + (x^5 + x^2)y = x^8 + x^7 + x^6 + 4x^5 + 3x^4 + 2x^3 + 4x^2 + 2x$$ of genus 4, with discriminant -1,064,000, which was found by Harrison (\cite{HarrisonList}). It has Mordell-Weil rank 0. We found $L(A,1) \approx 0.09889146$, $P_A \approx 178.0046$, $c_2 = 2$, $c_p = 1$ for all other $p$, and $|A(\Q)_{\mathrm{tors}}| = 60$, yielding $|\Sh|_{\mathrm{an}} = 1.0000000$. Again the torsion is not computed in a provable way. However by reducing modulo 3, we found that the torsion is a divisor of 180. As far as we are aware this is the first example of a curve of genus 4 for which BSD has been numerically verified.
\item The curve $$y^2 + (x^6 + x^4 + 1)y = x^4 + x^2$$ of genus 5, with discriminant 116,985,856, found in the aforementioned list. It has Mordell-Weil rank 0. We found $L(A,1) \approx 0.1002872$, $P_A \approx 579.2589$, $c_p = 1$ for all $p$, and $|A(\Q)_{\mathrm{tors}}| = 76$, yielding $|\Sh|_{\mathrm{an}} = 1.0000000$. As this curve does not have a rational Weiersta\ss\ point (which we actually do assume for most of the article), the search for torsion points was much more cumbersome, due to the Mumford representation not behaving well in this case. Again it is not provable; the best upper bound for the torsion that we found is 304. As far as we are aware this is the first example of a curve of genus 5 for which BSD has been numerically verified.
\end{itemize}

\begin{remark}
It could be the case that some of these curves have isomorphic (or isogenous) Jacobians. Then we actually verified BSD two times for the same abelian variety. In the verification process, we did not check for this.
\end{remark}

\begin{remark}
Even though for all our curves the verification went well, it should be remarked that problems are to be expected when trying to verify BSD for curves with higher discriminant (or rather, higher conductor). The computation of the $L$-function takes much longer in these cases. Also the computation of the regulator will be harder, as the heights of the points involved might increase, in particular in case the Mordell-Weil rank is higher.

It should be feasible to carry out the verification for more of the small examples from Harrison's list, \cite{HarrisonList} of genus 4, as long as the maximum bad prime is small enough. We also tried the verification for some more examples of genus 5, but in these cases the computation of the special value of the $L$-function was taking hours and the computation of the regular model sometimes did not seem to finish in reasonable time.
\end{remark}

\onlyinsubfile{}

\newcommand{\nf}{\Q}

\section{Theory of differentials}

Let $C/\nf$\check{Now this is over $\Q$. I could add a remark that everything generalises over number fields, except for the problem that some modules are not free, but only locally free.} be a smooth, geometrically irreducible, projective curve of genus $g$ over $\nf$. Let $J/\nf$ be its Jacobian. The goal of this section is to define the period of $J$, and to describe a way to compute it in the case $C$ is hyperelliptic. We will be following the algorithm described in \cite[sect.\ 3.5]{EmpiricalEvidence}.

First we will discuss both the theoretical considerations that are needed for this algorithm.

Throughout the section $p$ will be a prime %of the ring of integers $\mathcal{O}_K$ of $K$,
and $S$ will be the scheme $\Spec{(\Z_{(p)})}$. The generic point of $S$ is called $\eta$ and the special point $p$.

\subsection{Preliminaries}

First, for completeness, we will recall the following definition.

\begin{definition}[\textrm{\cite[p.\ 166]{BLR}}]
A {\em (relative) curve} $\mathcal{C}$ over $S$ is a normal, proper, flat $S$-scheme, such that for all $t \in S$, the scheme $\mathcal{C}_t = \mathcal{C} \times_S k(t)$ is of pure dimension 1. A {\em model of $C$ over $S$} is a relative curve $\mathcal{C}$ over $S$ together with an isomorphism $\mathcal{C}_{\eta} \cong C$.
\end{definition}

\begin{remark}
Without the normality assumption, the special fibre of a curve over $S$ could have embedded components. In order to be able to use the results from \cite{BLR}, which have been partially derived from \cite{Raynaud}, it is necessary to not have embedded components.
\end{remark}

Let $\mathcal{J}$ be a N\'eron model of $J$ over $S$, and let $\mathcal{C}/S$ be a regular model of $C$. Assume that the geometric multiplicities of the irreducible components of $\mathcal{C}_p$ in $\mathcal{C}_p$ have greatest common divisor 1.

\begin{theorem}[\textrm{\cite[Thm.\ 4(b), sect.\ 9.5, p.\ 267]{BLR}}]
Under these conditions, $\Pic_{\mathcal{C}/S}^0$ is a separated scheme and $\Pic_{\mathcal{C}/S}^0$ coincides with the identity component of $\mathcal{J}$.
\end{theorem}

From \cite[Prop.\ 5.2, p.\ 46]{Raynaud}, it now follows that $\mathcal{C}/S$ is cohomologically flat, which we will need for the next part.

										% Differentials %

\subsection{Differentials of Jacobian and regular model}

A classical theorem (see for example \cite[Prop.\ 2.2, p.\ 172]{MilneJacobian}) relates the differentials on the Jacobian of a smooth curve over a field with the differentials on the curve itself. We will generalise this to $\mathcal{J}$ and $\mathcal{C}$.

\begin{definition}[\textrm{\cite[Def.\ 4.7, sect.\ 6.4.2, p.\ 239]{Liu}}]\label{def:canonicalsheaf}
Let $Y/T$ be a quasi-projective locally noetherian scheme. Let $i : Y \rightarrow Z$ be an immersion into a smooth scheme $Z/T$. Then the {\em canonical sheaf of $Y/T$} is defined to be the $\O_Y$-module $$\omega_{Y/T} := \det(i^*(\mathcal{I}/\mathcal{I}^2))^{\vee} \otimes_{\O_T} i^*(\det \Omega^1_{Z/T}),$$
where $\mathcal{I}$ is the sheaf of ideals defining $Y$ in an open $Z' \subset Z$ containing $Y$ as closed subset. This is independent of the choice of $Z$ and $i$, see {\it loc.\ cit.}
\end{definition}

The following lemma generalises the aforementioned theorem.

\begin{lemma}\label{thm:DifferentialsCurveJacobian}
There are canonical isomorphisms of $\O_S$-modules 
$$\xymatrix{
\Omega^1_{\mathcal{J}/S}(\mathcal{J}) \ar[d]^{\sim} 	&&\omega_{\mathcal{C}/S}(\mathcal{C}) \ar[d]_{\sim}^{\textrm{GD}}	\\
\Hom_{\O_S}(\mathrm{Lie}(\mathcal{J}),\O_S) \ar[rr]^{\alpha}_{\sim}					&&\Hom_{\O_S}(R^1f_*(\O_{\mathcal{C}}), \O_S)	
}$$
\end{lemma}

\begin{proof}
The right hand isomorphism is given by Grothendieck duality, see \cite[Sect.\ 6.4.3, p.\ 243]{Liu}. The bottom isomorphism, $\alpha$, is from \textrm{\cite[Thm.\ 8.4.1, p.\ 231]{BLR}} (here we use that $\mathcal{C}/S$ is cohomologically flat). Getting the left hand isomorpism is a little bit more involved. 

First remark that global differentials on an abelian variety are invariant. As the image of $J$ is dense in $\mathcal{J}$, this also holds for the differentials in $\Omega^1_{\mathcal{J}/S}(\mathcal{J})$. Combining this with \cite[Prop.\ 4.2.1, p.\ 100]{BLR}, we get \begin{equation}\label{arrow:BLR100} \Omega_{\mathcal{J}/S}(\mathcal{J}) = \Omega_{\mathcal{J}/S}^1(\mathcal{J})^{\mathrm{inv}} =  e^*\Omega_{\mathcal{J}/S}^1(S),\end{equation} where $e : S \rightarrow \mathcal{J}$ is the unit section. Now, by \cite[Prop.\ 6.1.24, p.\ 217]{Liu}, we get an exact sequence of $\O_S$-modules
$$\mathfrak{m} / \mathfrak{m}^2 \rightarrow e^*\Omega_{\mathcal{J}/S}^1 \rightarrow \Omega_{S/S}^1 = 0,$$
where $\mathfrak{m}$ is the ideal of the schematic image of $e$ inside $\mathcal{J}$. As both $\mathfrak{m}/\mathfrak{m}^2$ and $\Omega^1_{\mathcal{J}/S}$, and hence $e^*\Omega^1_{\mathcal{J}/S}$ are locally free of rank $g$ (as $\mathcal{J}$ is regular), we get that the kernel of $\mathfrak{m}/\mathfrak{m}^2 \rightarrow e^*\Omega^1_{\mathcal{J}/S}$ is torsion. As $\mathcal{O}_S$ is torsion-free in our case, and hence the locally free module $\mathfrak{m}/\mathfrak{m}^2$ is torsion-free, we find a canonical isomorphism of $\O_S$-modules $$e^*\Omega^1_{\mathcal{J}/S} = \mathfrak{m}/\mathfrak{m}^2 = \ShHom_{\O_S}(\mathrm{Lie}(\mathcal{J}), \O_S),$$
which gives, by taking global sections and composing with eqn.\ \eqref{arrow:BLR100}, the construction of the left hand isomorphism in the diagram.
\end{proof}

\begin{remark}\label{rmk:identifications}
Under the natural identifications $\Omega^1_{\mathcal{J}/S}(\mathcal{J}) \otimes_{\Z_{(p)}} \Q = \Omega^1_{J/\Q}(J)$ and $\omega_{\mathcal{C}/S}(\mathcal{C}) \otimes_{\Z_{(p)}} \Q = \Omega^1_{C/\Q}(C)$, the isomorphism $\Omega^1_{\mathcal{J}/S}(\mathcal{J}) \cong \omega_{\mathcal{C}/S}(\mathcal{C})$ in the lemma above is compatible with the aforementioned classical isomorphism $\Omega^1_{J/\Q}(J) \cong \Omega^1_{C/\Q}(C)$.
\end{remark}

								% Algorithm period %

\subsection{Algorithm for the real period}

Suppose that $\omega_1, \ldots, \omega_g \in \Omega^1_{C/\Q}(C)$ are such that, for every prime $p$, they form a $\Z_{(p)}$-basis of $\omega_{\mathcal{C}/S}(\mathcal{C})$, under the identification $\omega_{\mathcal{C}/S}(\mathcal{C}) \otimes_{\Z_{(p)}} \Q = \Omega^1_{C/\Q}(C)$. In other words, cf.\ Lemma \ref{thm:DifferentialsCurveJacobian}, suppose that $\omega_1, \ldots, \omega_g$ correspond to generators of $\Omega^1_{\mathcal{J}_{\Z} / \Z}(\mathcal{J}_{\Z})$, where $\mathcal{J}_{\Z} / \Z$ is a N\'eron model of $J$ over $\Spec{\Z}$. Moreover, let $\gamma_1, \ldots, \gamma_{2g} \in H^1(C,\Z)$ form a symplectic basis for the homology. Then the real period can be defined as follows.

\begin{definition}
The {\em real period} of $C$ is the covolume of the lattice 
$$ \Z (a_1 + \overline{a_1}) + \ldots + \Z (a_{2g} + \overline{a_{2g}} ) \subset \R^g,$$
where $a_i = ( \int_{\gamma_i} \omega_j )_{j=1}^g \in \C^g$ for $i = 1, \ldots, 2g$.
\end{definition}

Now suppose that we are working with a hyperelliptic curve given by $y^2 = f$ for some $f \in \Q[x]$. Then, due to Van Wamelen there is a procedure in \texttt{Magma} to compute a symplectic basis of $H^1(C,\Z)$ as mentioned before, and the integrals $\int_{\gamma_i} \frac{x^{j-1} \cdot dx}{y}$ for all $i = 1, \ldots, 2g$ and $j = 1, \ldots, g$.

In order to compute the real period, we only need to find a basis $\omega_1, \ldots, \omega_g$ as above in terms of the differentials $\frac{x^{j-1} \cdot dx}{y}$. For our purpose, the calculation can be done for each prime $p$ separately.
Fortunately for us, due to Donnely, \texttt{Magma} also has an algorithm to compute explicit equations for a regular model $\mathcal{C}$ of $C$ over $S$. It will represent $\mathcal{C}/S$ by giving charts, each of which is a relative complete intersection. The following lemma explicitly gives the isomorphism $\omega_{\mathcal{C}/S}(\mathcal{C}) \otimes_{\Z_{(p)}} \Q \cong \Omega^1_{C/\Q}(C)$ that we need to compute whether a certain differential is vanishing or having a pole on one of the components of the special fibre (Step 5 and 6 in Algorithm \ref{algo}).

\newcommand{\Xeta}{\mathcal{X}_{\eta}}
\begin{lemma}\label{thm:ExplicitDifferentials}
Let $\mathcal{X} \subset \A^n_{S} = \Spec(\Z_p[x_1, \ldots , x_n])$ be regular, flat, and of relative dimension 1 over $S = \Spec{\Z_{(p)}}$. %\check{Is this enough?}
Suppose that $\mathcal{X}$ is a relative complete intersection inside $\A^n_S$, given by equations $g_1 = \ldots = g_{n-1} = 0$, with $g_i \in \Z_{(p)}[x_1, \ldots , x_n]$. Moreover, suppose that the generic fibre $\Xeta$ is smooth over $\Q$.

Then, on the one hand, after possibly reordering $x_1, \ldots, x_n$, we may and will assume that $\Omega^1_{k(\Xeta)/\Q}$ is a $k(\Xeta)$-vector space of dimension 1 generated by $dx_n$%\check{Add: up to reordering of $dx, dy, dz$.}
. This space contains $\Omega^1_{\Xeta/\Q}(\Xeta)$. On the other hand, we can define $\omega_{\mathcal{X}/S}$ using this immersion into $\A^n_S$ (cf.\ Def.\ \ref{def:canonicalsheaf}). 
Then $\omega_{\mathcal{X}/S}$ is free of rank 1 and generated by an element, which we will denote by ${(g_1 \wedge \ldots \wedge g_{n-1})}^{\vee} \otimes {dx_1 \wedge \ldots \wedge dx_n}$.
Then there is a canonical isomorphism of $\Q$-vector spaces
\[\Omega^1_{\Xeta/\Q}(\Xeta) \stackrel{\sim}\longrightarrow \omega_{\mathcal{X}/S}(\mathcal{X}) \otimes_{\Z_{(p)}} \Q,\]
which is given by
\[f \cdot dx_n \mapsto f \cdot  \det \begin{pmatrix}
\partial{g_i} / \partial{x_j}
\end{pmatrix}_{i,j=1}^{n-1}
 \cdot {(g_1 \wedge \ldots \wedge g_{n-1})}^{\vee} \otimes  {dx_1 \wedge \ldots \wedge dx_n}.\]
\end{lemma}

%\begin{proof}[Proof (original)]This lemma is just a matter of explicitly writing out some isomorphisms. See Lemma 6.4.5 (p.\ 238), Lemma 6.4.1 (p.\ 236), Cor.\ 6.3.22 (p.\ 233) Cor.\ 6.3.14 (p.\ 231) Prop.\ 6.3.13 (p.\ 230), Prop.\ 6.1.24 (p.\ 217), and Prop.\ 6.1.8 (p.\ 212) from \cite{Liu}.
%\end{proof}

\begin{proof}
On the one hand, we can consider $\Xeta \subset \Xeta$, on the other hand, we have an embedding $\Xeta \subset \A^n_{\Q}$. Both give us a way to construct $\Omega^1_{\Xeta/\Q}$, and \cite[Lem.\ 6.4.5, p.\ 238]{Liu} gives an explicit natural isomorphism between them. What is left to check, is that this isomorphism is exactly the one described in the statement of Lemma \ref{thm:ExplicitDifferentials}.

We will break down the proof of \cite[Lem. 6.4.5, p.\ 238]{Liu} to find the map explicitly. In this lemma, we will take $X = Z_1 = \Xeta$, $Y = \Spec{\Q}$ and $Z_2 = \A^n_{\Q}$, and we let $i_2\colon \Xeta \rightarrow \A^n_{\Q}$ be the map induced by the embedding of $\mathcal{X}$ into $\A^n_S$. The two exact sequences, induced by \cite[Cor.\ 6.3.22, p.\ 233]{Liu} are
$$0 \rightarrow 0 \rightarrow \mathcal{C}_{\Xeta/W} \rightarrow \mathop{i_2^*} \Omega^1_{\A^n_{\Q}/\Q} \rightarrow 0 \qquad \textrm{and} \qquad 0 \rightarrow \mathcal{C}_{\Xeta/\A^n_{\Q}} \rightarrow \mathcal{C}_{\Xeta/W} \rightarrow \Omega^1_{\Xeta/\Q} \rightarrow 0,$$
where $W = \Xeta \times_{\Q} \A^n_{\Q}$, and the map $h \colon \Xeta \rightarrow W$ is given by $(\mathrm{id}_{\Xeta},i_2)$, and $\mathcal{C}_{\Xeta/W} = \mathop{h^*} \mathcal{I}_h / \mathcal{I}_h^2$ and $\mathcal{C}_{\Xeta} / \A^n_{\Q} = \mathop{i_2^*} \mathcal{I}_{i_2} / \mathcal{I}_{i_2}^2$ with $\mathcal{I}_h$ and $\mathcal{I}_{i_2}$ the sheaf of ideals on $W$ and $\A_q^n$ respectively, defining $\Xeta$.

We will make the maps in these exact sequences explicit, starting with the first sequence. Let $p_1 \colon W \rightarrow \Xeta$ and $p_2 \colon W \rightarrow \A^n_{\Q}$ be the two projections. We know that $\Omega^1_{\A_{\Q}^n/\Q}$ is a free sheaf generated by $n$ elements $dx_1, \ldots, dx_n$. Now $\Omega^1_{W/\Xeta}$ is identified with $\mathop{p_2^*} \Omega^1_{\A^n_{\Q}/\Q}$, and in this identification the differential $dx_j$ is mapped to $dz_j$, where $z_j = \mathop{p_2^*} x_j$. By pulling back along $h$, we get an identification $\mathop{h^*} \Omega^1_{W/\Xeta} = \mathop{i_2^*} \Omega^1_{\A^n_{\Q}/\Q}$.

Now the isomorphism $\mathcal{C}_{\Xeta/W} \rightarrow \mathop{h^*} \Omega^1_{W/\Xeta}$ is ultimately coming from \cite[Prop.\ 6.1.8, p.\ 212]{Liu}. The sheaf $\mathcal{I}_h / \mathcal{I}_h^2$ has is generated by $z_j - y_j$, for $j = 1, \ldots, n-1$, where $y_j = \mathop{p_1^*} \mathop{i_2^*} x_j$. These are mapped to $d(z_j - y_j) = dz_j$ in $\mathop{h^*} \Omega^1_{W/\Xeta}$ or to $dx_j$ in $\mathop{i_2^*} \Omega^1_{\A_{\Q}^n/\Q}$.

To understand the morphism $\mathcal{C}_{\Xeta/\A^n_{\Q}} \rightarrow \mathcal{C}_{\Xeta/W}$ in the second sequence, we have to go back to \cite[Cor.\ 6.3.22]{Liu}. The sheaf $\mathcal{I}_{i_2} / \mathcal{I}_{i_2}^2$ is generated by the functions $g_1, \ldots, g_{n-1}$. Following the proof of the aforementioned corollary, we consider the following cartesian diagram.
$$\xymatrix{
\Xeta \times_{\Q} \Xeta \ar[rr]^{(\pi_2 , \mathop{\pi_1^*} i_2)} \ar[d]_{\pi_1} &&W \ar[d]^{p_2}	\\
\Xeta \ar[rr]^{i_2} &&\A^n_{\Q}
}$$
Here $\pi_1$ and $\pi_2$ are the first and second coordinate projections $\Xeta \times_{\Q} \Xeta \rightarrow \Xeta$. The map $h \colon \Xeta \rightarrow W$ from the bottom left to the top right, using the universal property of the product, gives rise to the diagonal section $\Delta \colon \Xeta \rightarrow \Xeta \times_{\Q} \Xeta$ of $\pi_1$. Then, there is the identification $$\mathcal{C}_{\Xeta/\A^n_{\Q}} = \mathop{\Delta^*} \mathop{\pi_1^*} \mathcal{C}_{\Xeta/\A^n_{\Q}} = \mathop{\Delta^*} \mathcal{C}_{\Xeta \times_{\Q} \Xeta / W},$$ identifying the functions $g_i$ in $\mathcal{C}_{\Xeta/\A^n_{\Q}}$ with the functions $\mathop{p_2^*} g_i$ in $\mathop{\Delta^*} \mathcal{C}_{\Xeta \times_{\Q} \Xeta/W}$. In other words, if you express the $g_i$ in terms of the variables $x_j$ on $\A^n_{\Q}$, then you get $\mathop{p_2^*} g_i$ by replacing all the $x_j$'s by $z_j$'s.

The map $\mathcal{C}_{\Xeta/W} \rightarrow \Omega^1_{\Xeta/\Q}$ is constructed in an analogous way to the construction of the map $\mathcal{C}_{\Xeta/W} \rightarrow \mathop{i_2^*} \Omega^1_{\A^n_{\Q}/\Q}$. It sends $z_j - y_j$ to $-dw_j$, where $w_j = \mathop{i_2^*} x_j$, on $\Omega^1_{\Xeta/\Q}$. 

Now the isomorphism 
\begin{align*}
&\det \mathcal{C}_{\Xeta/\A^n_{\Q}} \otimes \det \Omega^1_{\Xeta/\Q} &&\longrightarrow \qquad \det \mathcal{C}_{\Xeta/W} &&\longrightarrow \qquad \det \mathop{i_2^*} \Omega^1_{\A^n_{\Q}/\Q}\\
&(g_1 \wedge \ldots \wedge g_{n-1}) \otimes dw_n &&\longmapsto p_2^*g_1 \wedge \ldots \wedge p_2^* g_{n-1} \wedge dw_n &&\longmapsto dg_1 \wedge \ldots \wedge dg_{n-1} \wedge dx_n
\end{align*}
is constructed cf.\ \cite[Lem.\ 6.4.1, p.\ 236--237]{Liu}. Of course,
$$dg_1 \wedge \ldots \wedge dg_n \wedge dx_n = \det \begin{pmatrix} \partial g_i / \partial x_j \end{pmatrix}_{i,j=1}^{n-1} \cdot dx_1 \wedge \ldots \wedge dx_n.$$
Recall that $\omega_{\mathcal{X}/S} = \det(\mathop{\iota^*} \mathcal{I}_{\iota} / \mathcal{I}_{\iota}^2)^{\vee} \otimes_{S} \, \mathop{\iota^*} \det \Omega^1_{\mathcal{X}/S}$, where $\iota \colon \mathcal{X} \rightarrow \A^n_{S}$ is the embedding, and $\mathcal{I}_{\iota}$ is the sheaf of ideals on $\O_{\A_S^n}$ defining $\mathcal{X}$. After base change to $\Q$, this becomes $(\det \mathcal{C}_{\Xeta/\A^n_{\Q}})^{\vee} \otimes_{\Q} \det \mathop{i_2^*} \Omega^1_{\A^n_{\Q}/\Q}$.
The result now follows immediately.
\end{proof}

Altogether, this leads to the following algorithm.

\begin{algorithm}\mbox{}\\[-0.9cm]\label{algo}
\begin{itemize}
\item[] {\em Input}: monic polynomial $f \in \Z[X]$ of degree $2g+1$ describing a hyperelliptic curve $C$ of genus $g$ over $\Q$.
\item[] {\em Output}: the period $\Omega$ of $C$.
\item[] {\em Step 1}: calculate the so-called big period matrix $(\int_{\gamma_i} \omega_j)_{i=1, \ldots, 2g, \, j=1, \ldots g}$ of $J$, where the notation is as before, using the \texttt{Magma} command \verb+BigPeriodMatrix+ (due to Van Wamelen).
\item[] {\em Step 2}: for each subset $I \subset \{1, \ldots, 2g\}$ with $|I| = g$, calculate the covolume $P_I := \left|\det\left(\int_{\gamma_i} \omega_j + \overline{\int_{\gamma_i} \omega_j}\right)_{i \in I,\, j=1,\ldots, g}\right|$.
\item[] {\em Step 3}: use Euclid's algorithm to find a generator $P$ for the lattice spanned by the $P_I$.
\item[] {\em Step 4}: for each bad prime $p$, calculate a regular model $\mathcal{C} / \Z_{(p)}$ of $C$, using the \texttt{Magma} command \verb+RegularModel+. This will give us a representation of $\mathcal{C}$ by charts which are relative complete intersections.
\item[] {\em Step 5}: for each of the differentials $\omega_1, \ldots, \omega_g$, check that if it has a pole on any of the irreducible components of the special fibre of $\mathcal{C}$. If so, adjust the basis by multiplying the differential having a pole with $p$ to get a new basis $\underline{\omega}'$ and apply Step 5 again (until the basis is not changing anymore).
\item[] {\em Step 6}: for each $(c_j)_{j=1}^g \in \{0, \ldots, p-1\}^g \setminus \{(0,0, \ldots, 0)\}$, check if $\sum_j c_j \omega_j$ vanishes on the whole special fibre of $\mathcal{C}$. If so, adjust the basis $\underline{\omega}'$ by replacing one of the $\omega_j$ such that $c_j \neq 0$ with $\tfrac1p \sum_j c_j \omega_j$, then apply Step 6 again (until the basis is not changing anymore).
\item[] {\em Step 7}: for each bad prime $p$ compute $p^{a-b}$, where $a$ is the number of basis adjustments done in Step 5, and $b$ is the number of basis adjustments done in Step 6 (this is also the determinant of the change of basis matrix whose columns express $\underline{\omega}'$ in terms of $\underline{\omega}$). Then take the product $W$ over $p$ of these determinants, and output $W \cdot P$.
\item[] {\em End.}
\end{itemize}
\end{algorithm}

\onlyinsubfile{}

\section{Computation of other terms in BSD formula}

Throughout this section, we will use the following notation.

\begin{notation}
We define $H/\Q$ to be a hyperelliptic curve of genus $g$. When a prime $p$ is introduced, $\mathcal{H} / \Z_{(p)}$ is a regular model of $H$ over $\Z_{(p)}$. The Jacobian of $H$ is denoted by $J$, and the N\'eron model of $J$ over $\Z$ is called $\mathcal{J}$.
\end{notation}

Moreover, we will assume that $H$ is given by a model of the form $y^2 = f(x)$, where the input polynomial $f(x)$ has odd degree (and hence $H$ has a rational Weierstra\ss\ point).

\onlyinsubfile{\section{Computational methods}}

\subsection{Torsion subgroup and rank}				% TORSION SUBGROUP %

In order to compute the torsion group and algebraic rank of $J$, we will be computing upper and lower bounds.

For the torsion, upper bounds are given by considering the reduction of $J$ at good primes. For the algebraic rank, upper bounds are given by considering 2-Selmer groups. This is already implemented in \texttt{Magma} by Stoll.

To get lower bounds, we try to find as many points as possible on $J$. For genus 2, this is already implemented in \texttt{Magma}. For genus 3, 4 and 5, the author implemented a simple search algorithm for points, using the Mumford representation that \texttt{Magma} is using to represent points on $J$.

In fact, for Jacobians of curves $J$ and $J^{\vee}$ are isomorphic. Hence, in order to verify the BSD conjecture up to squares in this case, it is actually not necessary to know the size of the torsion subgroup at all.

%\begin{example}
%The torsion subgroup of $J_1$ is isomorphic to $\Z/5\Z$.
%\end{example}

%\begin{example}
%The torsion subgroup of $J_2$ is isomorphic to $\Z/2\Z$.
%\end{example}

%\begin{example}
%In our search for points of small height on $J_3$ we found three points of order 2 and 1 point of order 1. Therefore, for further calculations we will assume that the torsion subgroup is isomorphic to $\Z/2\Z \times \Z/2\Z$.\check{Maybe do reductions modulo primes to get an upper bound on the torsion. I can get to 16 now.}
%\end{example}

\onlyinsubfile{}

\onlyinsubfile{\section{Computational methods}}

\subsection{$L$-function}\label{sec:Lfunc}					% L-FUNCTION %

In this section, we will briefly discuss the computation of the special value of the $L$-function associated to the Jacobian of a hyperelliptic curve. For a complete definition and theoretical background on the $L$-function, see \cite{Serre}.

The idea used to compute the $L$-function is as follows. The local $L$-factors at the good primes $p > 2$ can be found by counting points in $\mathcal{J}_p(\F_{p^m})$ for sufficiently many $m \geq 1$. In order to find the local $L$-factors at the bad places, one uses the functional equation. The idea is to guess, in a clever way, the conductor and, for the bad primes, the local $L$-factors, in such a way that the $L$-function obtained satisfies the conjectural functional equation, see also \cite[sect.\ 5, p.\ 243--245]{Database2}.

\begin{remark}\label{rmk:Ogg}
To guess the 2-part of the conductor, the following naive version of Ogg's formula is used:
$$f^{\textrm{guess}} = v(\Delta) - n + 1.$$
Here, $v(\Delta)$ is the valuation of the (naive) minimal discriminant, $n$ is the number of geometrically irreducible components in a minimal regular model, and $f^{\textrm{guess}}$ is our guess for the 2-valuation of the conductor. The formula, in this shape, does not give the correct 2-valuation of the conductor in general. For curves of genus 2 over a henselian discrete valuation ring with algebraically closed residue field, we can deduce the formula
$$f = v(\Delta) - n + 1 - 11 \cdot c(X),$$
from \cite{LiuConductor}, where $c(X)$, as defined in loc.\ cit., is a non-negative integer. Over general discrete valuation rings, the discriminant could change after a quadratic field extension, cf.\ \cite[Prop.\ 4, p.\ 4595]{LiuModels}. In this case, it drops by $2(2g+1)$. So, for genus 2, in case $v(\Delta) < 10$, the discriminant will apparently not change anymore, and $c(X) = 0$ must hold for the 2-valuation $f$ of the conductor to not become negative. Hence, the naive version of Ogg's formula holds in this case.
\end{remark}

In \cite{ComputeL}, Tim Dokchitser describes a trick with an inverse Mellin transform in order to actually evaluate the $L$-function. This has been implemented by him, together with Vladimir Dokchitser, in \texttt{Magma}. This is the method we used for our calculations. However, it is useful to remark that the runtime increases quickly when the conductor increases and that this could probably by remedied by using the methods from \cite{HMS}.

%\begin{example}
%For $H_1$ we found that $$L(J_1, 1) \approx 2.08419385369113888173282768910.$$ This agrees with the fact the algebraic rank of $J_1$ is 0.
%\end{example}

%\begin{example}
%For $H_2$ we found that $$L(J_2, 1) \approx 5 \cdot 10^{-31} \textrm{ and } L'(J_2,1) \approx 0.819558937768934171069200441694.$$
%This agrees with the fact that the algebraic rank of $J_2$ is 1. 
%\end{example}

%\begin{example}
%For $H_3$ we found that $$L(J_3, 1) \approx 3 \cdot 10^{-30} \textrm{ and } L'(J_3, 1) \approx 0.869490085404871842347716239330.$$
%This agrees with the fact that the algebraic rank of $J_3$ is 1. {\color{red} Maybe some remark about the forcing of Ogg's formula at 2.}
%\end{example}

\onlyinsubfile{}

%\subfile{Period}

\onlyinsubfile{\section{Computational methods}}
\subsection{Regulator}						% REGULATOR %

Using the points on $J$ that we found when computing the algebraic rank, we will compute the regulator. In order to do that, we need to calculate the height pairing for several pairs of points.

Due to work of Holmes (\cite{HeightHolmes}) and M\"uller (\cite{HeightMueller}) it is now known how arithmetic intersection theory could be used to do this calculation. This has also been implemented in \texttt{Magma} for hyperelliptic curves by M\"uller, and works in practice for genus up to 10.

In many cases, especially in genus 3, 4 and 5, the height bound we use for point finding is not high enough to provably compute the regulator. The upper bounds for difference between the naive and canonical height are quite big in some cases, see for example \cite{MullerStoll} for genus 2.
 In that case, we can only obtain a finite index subgroup of the Mordell-Weil group. Therefore, the regulator that we get might be a square multiple of the actual regulator of $J$. Hence, the conjectural order of $\Sh$, assuming BSD, might be a multiple of the order that we compute. 

%\begin{example}
%For $H_1$, Magma tells us that $J_1$ has algebraic rank 0. Hence, the regulator is 1.
%\end{example}

%\begin{example}
%For $H_2$, Magma tells us that $J_2$ has algebraic rank 1. Having looked at points of height up to {\color{red} (...)}, we found that the point $(x-1, 1)$, in Mumford representation, is probably a generator for the free part. The heigth pairing of the point with itself is $$0.0474425704742192075988905184458.$$
%The other values for the height pairing of points with itself are all square multiples of this number, which indicates that this point is very likely a generator for the free part. Hence, the regulator probably equals the aforementioned number.
%\end{example}

%\begin{example}
%For $H_3$, Magma tells us that $J_3$ has algebraic rank 1. Having looked at points of height up to {\color{red} (...) pretty small}, we found that the point $(x,-1)$, in Mumford representation, is probably a generator for the free part. The height pairing of the point with itself is $$0.230972622186586831975393030694.$$
%The other values for the height pairing of points with itself are all square multiples of this number\check{There is one exception, which is very strange. Is the height pairing supposed to be bilinear? Now I get this: HeightPairing(-4*P11, -4*P11) =
%4.042...
%HeightPairing(P11,P11) = 0.230...
%}, which indicates that this point is very likely a generator for the free part. Hence, the regulator probably equals the aforementioned number.
%\end{example}

\onlyinsubfile{}

\renewcommand{\ker}{\mathop{\mathrm{Ker}}}
\onlyinsubfile{\section{Computational methods}}

\subsection{Tamagawa numbers}				% TAMAGAWA NUMBERS %

Suppose that we have a regular model $\mathcal{H}^{\mathrm{s}}$ of $H$ over the strict henselisation of $\Z_{(p)}$. Then in \cite[Thm.\ 1.1, p.\ 277]{BoschLiu}, Bosch and Liu give an exact sequence
$$0 \rightarrow \mathop\mathrm{Im} \overline{\alpha} \rightarrow \ker \overline{\beta} \rightarrow \phi_A(\overline{\F_p}) \rightarrow 0$$
of $\mathrm{Gal}(\overline{\F_p} / \F_p)$-modules. Here $\phi_A(\overline{\F_p})$ is the geometric component group of the N\'eron model of $\mathcal{J}$. The map $\overline{\alpha} \colon \Z^{\overline{I}} \rightarrow \Z^{\overline{I}}$, with $\overline{I}$ indexing the components $\{\Gamma_i : i \in I\}$ of the special fibre of $\mathcal{H}^{\mathrm{s}}$, maps each component $\Gamma_j$ to $\sum_{i \in \overline{I}} e_i^{-1} \langle \Gamma_j, \Gamma_i \rangle \cdot \Gamma_i$, where $\langle \cdot, \cdot \rangle$ is the intersection pairing and $e_i$ is the geometric multiplicity of $\Gamma_i$ (in itself, which is 1 in our case). The map $\overline{\beta} \colon \Z^{\overline{I}} \rightarrow \Z$ maps each component $\Gamma_j$ to $d_j e_j$, where $d_j$ is the multiplicity of $\Gamma_j$ in the special fibre. Here, the Galois group $\mathrm{Gal}(\overline{\F_p}/\F_p)$ acts on $\Z^{\overline{I}}$ by its natural action on the components of the special fibre.

Due to Donnely, \texttt{Magma} is able to compute this geometric component group using this theorem, and moreover, because explicit equations exist for a regular model $\mathcal{H}$ of $H$ over $\mathbb{Z}_{(p)}$, we are able to compute the action of Frobenius on $\mathop\mathrm{Im} \overline{\alpha}$ and $\ker \overline{\beta}$.

The way regular models are constructed in \texttt{Magma} is by repeatedly blowing up non-regular points until the fibred surface is regular. To compute the Galois action on the components of the special fibre, we traced down this blow-up procedure, and in each step we computed the action of Galois on the points blown-up, and on the new components which appeared in the special fibre on the new blown-up charts.

The result is an implementation of a \texttt{Magma} package on top of the existing regular models package, which computes the action of the Galois group on $\phi_A(\overline{\F_p})$, and then computes the Tamagawa number, the order of $\phi_A(\F_p)$. The source code for this package will be released together with this article. It has been used to compute Tamagawa numbers for almost all of the 66,158 genus 2 curves present in \cite{lmfdb} (see also \cite{Database2}). This computation was finished within a few hours.

%\begin{example}
%The curve $H_1$ has only one bad prime, which is 5. The Tamagawa number at 5 is 5.\check{Maybe some detailed explanation?}
%\end{example}

%\begin{example}
%The curve $H_2$ has three bad primes: 2, 5 and 7. The Tamagawa number at 2 is 4, the Tamagawa numbers at 5 and 7 are 1.
%\end{example}

%\begin{example}
%The curve $H_3$ has three bad primes: 2, 3 and 31. The Tamagawa number at 2 is 2, the Tamagawa numbers at 3 and 31 are 1. 
%\end{example}

\renewcommand{\ker}{\mathop{\mathrm{ker}}}
\onlyinsubfile{}

\onlyinsubfile{\section{Computational methods}}

\subsection{Tate-Shafarevich group}		% TATE-SHAFAREVICH GROUP %

For our calculations, we do not calculate the order of the Tate-Shafarevich group. Instead, we only check whether the conjectural order, given by the BSD conjecture, is (up to a certain precision) a rational square or two times a rational square (with a small denominator) according to the criteria described in \cite{PoonenStoll}.

%\begin{example}
%For $H_1$, the BSD formule would give
%$$|\Sh(J_1)| = \frac{L^{(r)}(J_1,1) \cdot |J_1(\Q)_{\tors}|^2}{r! \cdot P_{J_1}R_{J_1} \cdot \prod_p c_p} \approx 4.00000000000000000000000000001,$$
%which is almost equal to 4, which agrees with the results from \cite{EmpiricalEvidence}.
%\end{example}

%\begin{example}
%For $H_2$, the BSD formule would give
%$$|\Sh(J_2)| = \frac{L^{(r)}(J_2,1) \cdot |J_2(\Q)_{\tors}|^2}{r! \cdot P_{J_2}R_{J_2} \cdot \prod_p c_p} \approx 1.00000000000000000000000000049,$$
%which is almost equal to 1, which agrees with the fact that Magma tells us that the order of $\Sh(J_2)$ is expected to be a square.
%\end{example}

%\begin{example}
%For $H_3$, the BSD formule would give
%$$|\Sh(J_3)| = \frac{L^{(r)}(J_3,1) \cdot |J_3(\Q)_{\tors}|^2}{r! \cdot P_{J_3}R_{J_3} \cdot \prod_p c_p} \approx 1.00000000000000000000000000003,$$
%which is almost equal to 1, which agrees with the fact that Magma tells us that the order of $\Sh(J_2)$ is expected to be a square.
%\end{example}

\onlyinsubfile{}

\end{document}